\documentclass[12pt]{article}

\usepackage[latin1]{inputenc}

\usepackage{pgf,tikz}
\usetikzlibrary{arrows}
\usetikzlibrary{patterns}

\usepackage{subfigure}

\usepackage{amssymb}

\usepackage{amsmath,epsfig}
\def\WE{\textrm{WE}}
\def\NS{\textrm{NS}}
\def\NE{\textrm{NE}}
\def\SE{\textrm{SE}}
\def\NW{\textrm{NW}}
\def\SW{\textrm{SW}}
\def\XY{\textrm{XY}}

\newtheorem{Theorem}{Theorem}
\newtheorem{Corollary}[Theorem]{Corollary}
\newtheorem{Proposition}[Theorem]{Proposition}

\newtheorem{Lemma}[Theorem]{Lemma}

\newcommand{\beq}{\begin{equation}}
\newcommand{\eeq}{\end{equation}}

\def\qed{\hfill$\Box$}

\usepackage{multirow}

\def\mybox{\hbox to 12.0pt}

\def\mybigbox{\hbox to 35.0pt}
\def\myverybigbox{\hbox to 60.0pt}

\def\ds{\displaystyle}    
\def\ov{\overline}
\def\sc{\scriptstyle}

\def\row{{\rm row}}
\def\col{{\rm col}}
\def\wgt{{\rm wgt}}

\def\xwgt{{\rm xwgt}}
\def\qxwgt{{\rm qxwgt}}

\def\x{{\bf x}}
\def\y{{\bf y}}

\def\1{{\bf 1}}

\def\+{\!\!+\!\!}

\def\magenta{\textcolor{magenta}}

\newcommand{\YT}[3]{
\vcenter{\hbox{
\begin{tikzpicture}[x={(0in,-#1)},y={(#1,0in)}] 
\foreach \k [count=\i] in {#3} {
 \foreach \e [count=\j] in \k {
  \draw (\i,\j) rectangle +(-1,-1);
  \draw (\i-0.5,\j-0.5) node {$#2\e$};
 }
}
\end{tikzpicture}
}}
}

\newcommand{\SYT}[3]{
\vcenter{\hbox{
\begin{tikzpicture}[x={(0in,-#1)},y={(#1,0in)}] 
\foreach \k [count=\i] in {#3} {
 \foreach \e [count=\j] in \k {
  \draw (\i,\j+\i-1) rectangle +(-1,-1);
  \draw (\i-0.5,\j+\i-1-0.5) node {$#2\e$};
 }
}
\end{tikzpicture}
}}
}

\newcommand{\wideSYT}[4]{
\vcenter{\hbox{
\begin{tikzpicture}[x={(0in,-#1)},y={(#2,0in)}] 
\foreach \k [count=\i] in {#4} {
 \foreach \e [count=\j] in \k {
  \draw (\i,\j+\i-1) rectangle +(-1,-1);
  \draw (\i-0.5,\j+\i-1-0.5) node {$#3\e$};
 }
}
\end{tikzpicture}
}}
}

\newcommand{\ASM}[3]{
\vcenter{\hbox{
\begin{tikzpicture}[x={(0in,-#1)},y={(#1,0in)}] 
\foreach \k [count=\i] in {#3} {
 \foreach \e [count=\j] in \k {
  \draw (\i-0.5,\j-0.5) node {$#2\e$};
 }
}
\end{tikzpicture}
}}
}

\newcommand{\CPM}[4]{
\vcenter{\hbox{
\begin{tikzpicture}[x={(0in,-#1)},y={(#2,0in)}] 
\foreach \k [count=\i] in {#4} {
 \foreach \e [count=\j] in \k {
  \draw (\i-0.5,\j-0.5) node {$#3\e$};
 }
}
\end{tikzpicture}
}}
}

\title{A note on corollaries to Tokuyama's Identity for symplectic Schur $Q$-Functions}

\author{
Ang\`ele M. Hamel\thanks{ 
Department of Physics and Computer Science,
Wilfrid Laurier University,
Waterloo, Ontario, N2L 3C5, Canada ({\tt ahamel@wlu.ca})}
\and 
Ronald C. King\thanks{
Mathematical Sciences, University of Southampton, 
Southampton SO17 1BJ, England ({\tt R.C.King@soton.ac.uk})}}

\begin{document}

\maketitle

\begin{abstract}
We present some corollaries to a symplectic primed shifted tableaux version of Tokuyama's identity expressed in terms
of other combinatorial constructs, namely generalised $U$-turn alternating sign matrices and strict symplectic Gelfand-Tsetlin patterns.
\end{abstract}

\section{Introduction}
\label{sec:intro}

Stimulated by a query from Vineet Gupta, Uma Roy, and Roger Van Peski we present a note on various ways of expressing the type $C$ Tokuyama identity 
that was first established in its most general form in the context of a symplectic primed shifted tableaux model~\cite{HK-bij}. This form,
Proposition 1.2 of~\cite{HK-bij}, involves parameters $\x=(x_1,x_2,\ldots,x_n)$, $\y=(y_1,y_2,\ldots,y_n)$ and $t$. Of these the
parameter $t$ carries very little significance, and may be eliminated as shown below in Corollary~\ref{cor-Tok-Q}. This result
can then be restated in terms of unprimed symplectic shifted tableaux to give Theorem~\ref{thm-main}. 
We then show in Section~\ref{sec-Tok-A-GT} that this theorem may be rewritten in terms of $U$-turn alternating sign matrices as in Corollary~\ref{cor-A} 
and in terms of strict symplectic Gelfand-Tsetlin patterns as in Corollary~\ref{cor-TokGT}, where the connection between symplectic shifted tableaaux 
and both $U$-turn alternating sign matrices and strict symplectic Gelfand-Tsetlin patterns is spelled out in Section~\ref{sec-A-GT}.
Finally the specialisation $y_k=qx_k$ for $k=1,2,\ldots,n$ is made in section~\ref{sec-yqx}, leading to four further type $C$ Tokuyama identities, 
including one due to Gupta {\em et al.}~\cite{GRVP} given here as our final Corollary~\ref{cor-TokGTq-final}.

This is intended as a brief note for those already familiar with the ideas and notations. However, for background information and motivation on Tokuyama's identity, see, for example, Tokuyama's original paper \cite{Tokuyama}, Okada's proof and variations \cite{Okadapartial, Okada}, and related work of Hamel and King \cite{HKWeyl, HK-bij}, Brubaker, Bump, and Friedberg \cite{Brubaker1}, Gupta, Roy, and Van Peski \cite{Gupta}, and Brubaker and Schultz \cite{BrubakerSchultz}, as well as the references therein.

\section{Background}
\label{sec:background}

For a partition $\mu$ of length $\ell(\mu)\leq n$ let ${\cal T}^\mu(n)$ be the set of symplectic tableaux $T$ of shape $F^\mu$
defined with respect to the alphabet 
\begin{equation}
    I=\{1<\ov{1}<2<\ov{2}<\cdots<n<\ov{n}\}\,,
\end{equation}
with the entries of $T$ taken from $I$:\\ 
\begin{tabular}{rl}
T1&weakly increasing across each row from left to right;\\
T2&strictly increasing down each column from top to bottom;\\
T3&such that $k$ or $\ov{k}$ appear no lower than the $k$th row.\\
\end{tabular} \\
This is 
exemplified in the case $n=4$ and $\mu=(4,3,3)$ by
\begin{equation}
\label{eqn-T}
T\ =\ \YT{0.2in}{}{
 {{1},{\ov1},{2},{\ov4}},
 {{\ov3},{4},{4}},
 {{4},{\ov4},{\ov4}},
}
\end{equation}
The deformed symplectic character $sp_\mu(\x;t)$ is defined by
\begin{equation}\label{eqn-wgtT}
   sp_\mu(\x;t)= \sum_{T\in{\cal T}^\mu(n)} \ \wgt(T)
\end{equation}
where $\wgt(T)$ is the product of weights $x_k$ and $t^2\ov{x_k}=t^2/x_k$ for each
entry $k$ and $\ov{k}$, respectively, in $T$.

Then for each strict partition $\lambda$ of length $\ell(\lambda)=n$ let ${\cal ST}^\lambda(n)$ be the set of 
symplectic shifted tableaux $ST$ of shifted shape $SF^\lambda$ defined with respect to the same alphabet $I$
with the entries of $ST$ taken from $I$:\\  
\begin{tabular}{rl}
ST1&weakly increasing across each row from left to right;\\
ST2&weakly increasing down each column from top to bottom;\\
ST3&strictly increasing down each diagonal from top left to bottom right.\\
\end{tabular} \\
This is 
exemplified below on the left in the case $n=5$ and $\lambda=(9,7,6,2,1)$.

Such tableaux may be refined through the addition of marks or primes.  
Let ${\cal QT}^\lambda(n)$ be the set of primed symplectic shifted tableaux $QT$ 
obtained from all $ST\in{\cal ST}^\lambda(n)$ by adding primes to entries in all possible ways such that\\
\begin{tabular}{rl}
QT1&any entry identical to its immediate left-hand neighbour must be unprimed;\\
QT2&any entry immediately above an identical entry must be primed;\\
QT3&any entry with neither of the above neighbours may be either unprimed or primed.\\ 
\end{tabular}

For the given tableau $ST$ shown below on the left, one such corresponding primed tableau $QT$ is
shown on the right.
\begin{equation}
\label{eqn-ST}
ST\ = \ 
\SYT{0.2in}{}{
 {{1},{\ov1},{2},{\ov2},{3},{3},{4},{\ov4},{\ov5}},
 {{2},{2},{\ov2},{\ov3},{4},{4},{\ov4}},
 {{\ov3},{4},{\ov4},{\ov4},{\ov4},{\ov4}},
 {{\ov4},{\ov4}},
 {{5}}
}
\qquad
QT\ = \
\SYT{0.2in}{}{
 {{1},{\ov1},{2'},{\ov2'},{3},{3},{4'},{\ov4'},{\ov5}},
 {{2'},{2},{\ov2},{\ov3'},{4},{4},{\ov4'}},
 {{\ov3},{4},{\ov4'},{\ov4},{\ov4},{\ov4}},
 {{\ov4},{\ov4}},
 {{5}}
}
\end{equation}

A deformed symplectic $Q$-function $Q_\lambda(\x/\y;t)$ is then defined by
\begin{equation}\label{eqn-wgtQ}
   Q_\lambda(\x/\y;t)= \sum_{QT\in{\cal QT}^\lambda(n)} \ \wgt(QT)
\end{equation}
where $\wgt(QT)$ is the product of the weights $x_k$, $y_k$, $t^2x_k^{-1}$ 
and $t^2y_k^{-1}$ for each entry $k$, $k'$, $\ov{k}$ and $\ov{k}'$, respectively, in $QT$.

Our starting point is Proposition~1.2 of~\cite{HK-bij} which can be written in the form
\begin{Proposition}
\label{prop-Tok-Qt}
Let $\mu$ be a partition of length $\ell(\mu)\leq n$ and $\delta= (n, n-1, \ldots, 1)$, so that
$\lambda=\mu+\delta$ is a strict partition of length $\ell(\lambda)=n$. 
Then for all $\x=(x_1,x_2,\ldots,x_n)$, $\y=(y_1,y_2,\ldots,y_n)$ and $t$
\begin{equation}\label{eqn-Qt-lambda}
Q_{\lambda}(\x/\y;t)= Q_{\delta}(\x/\y;t)\ sp_\mu(\x;t)\,,
\end{equation}
with
\begin{equation}\label{eqn-Qt-delta}
    Q_{\delta}(\x/\y;t) = \prod_{1\leq i\leq j\leq n} \ (x_i+y_j)(1+t^2x_i^{-1}y_j^{-1})\,.
\end{equation}
\end{Proposition}

The parameter $t$ is to some extent redundant. If for all $k=1,2,\ldots,n$ one maps $x_k$ and $y_k$ to $tx_k$ and $ty_k$, respectively, 
so that $t^2/x_k$ and $t^2/y_k$ map to $t/x_k$ and $t/y_k$, and takes into account the fact that every term in both $\wgt(T)$ and 
$\wgt(Q)$ of (\ref{eqn-wgtT}) and (\ref{eqn-wgtQ}) then carries a factor of $t$ to the power $|\mu|$ and $|\lambda|$, respectively,
it follows that in Proposition~\ref{prop-Tok-Qt} we have:
\begin{equation}
Q_{\lambda}(t\x/t\y;t)= t^{|\lambda|}\,Q_{\lambda}(\x/\y)\;~~Q_{\delta}(t\x/t\y;t)= t^{|\delta|}\,Q_{\delta}(\x/\y)\,;~~
           t^{|\mu|}\,sp_\mu(\x;1)=sp_\mu(\x)\,,
\end{equation}
where $Q_{\lambda}(\x/\y):=Q_{\lambda}(\x/\y;1)$, $Q_{\delta}(\x/\y):=Q_{\delta}(\x/\y;1)$ and $sp_\mu(\x):=sp_\mu(\x;1)$.
Then, since $|\lambda|=|\mu|+|\delta|$, we have
\begin{Corollary}
\label{cor-Tok-Q}
Let $\mu$ be a partition of length $\ell(\mu)\leq n$ and $\delta= (n, n-1, \ldots, 1)$, so that
$\lambda=\mu+\delta$ is a strict partition of length $\ell(\lambda)=n$.
Then for all $\x=(x_1,x_2,\ldots,x_n)$ and $\y=(y_1,y_2,\ldots,y_n)$
\begin{equation}\label{eqn-Q-lambda}
Q_{\lambda}(\x/\y)= \prod_{1\leq i\leq j\leq n} \ (x_i+y_j)(1+x_i^{-1}y_j^{-1})\ sp_\mu(\x)\,,
\end{equation}
\end{Corollary}
This $t=1$ version of Proposition~\ref{prop-Tok-Qt} has the advantage of being a Tokuyama-type identity for true, undeformed symplectic characters.

Although both versions have been expressed in terms of primed symplectic shifted tableaux it is a simple matter to write down the 
equivalent results expressed in terms of unprimed symplectic shifted tableaux. In particular the $t=1$ version takes the form:

\begin{Theorem}\label{thm-main}
Let $\mu$ be a partition of length $\ell(\mu)\leq n$ and $\delta= (n, n-1, \ldots, 1)$, so that
$\lambda=\mu+\delta$ is a strict partition of length $\ell(\lambda)=n$. 
Then for all $\x=(x_1,x_2,\ldots,x_n)$ and $\y=(y_1,y_2,\ldots,y_n)$
\begin{equation}\label{eqn-Tok-ST}
    ST_{\lambda}(\x/\y)= \prod_{1\leq i\leq j\leq n}\!\!\!(x_i+y_j)(1+x_i^{-1}y_j^{-1})\ \ sp_\mu(\x)\,,
\end{equation}
where
\begin{equation}\label{eqn-wgtPT}
   ST_{\lambda}(\x/\y) = \sum_{ST\in{\cal ST}^\lambda(\x)}\ \prod_{(i,j)\in{SF^\lambda}}\!\!\!\wgt(s_{ij}) 
	\quad\hbox{and}\quad 
	sp_\mu(\x)= \sum_{T\in{\cal T}^\mu(\x)}\ \prod_{(i,j)\in{F^\mu}}\!\!\!\wgt(t_{ij}) \,,
\end{equation}
with the weights of the entries $s_{ij}$ and $t_{ij}$ at position $(i,j)$ in $ST$ and $T$, respectively,
given by
\begin{equation}\label{eqn-swgt-twgt}
\begin{array}{|c|c|l|}
\hline
s_{ij}&\wgt(s_{ij})&\cr
\hline
 k&x_k &\hbox{if $s_{i,j-1}=k$}\cr
  &y_k &\hbox{if $s_{i+1,j}=k$}\cr
	&x_k+y_k&\hbox{otherwise}\cr
\hline
 \ov{k}^{\phantom{1}}&\ov{x}_k &\hbox{if  $s_{i,j-1}=\ov{k}$}\cr
  &\ov{y}_k &\hbox{if $s_{i+1,j}=\ov{k}$}\cr
	&\ov{x}_k+\ov{y}_k&\hbox{otherwise}\cr
\hline
\end{array}
\qquad\hbox{and}\qquad
\begin{array}{|l|l|}
\hline
t_{ij}&\wgt(t_{ij})\cr
\hline
 k&x_k\cr 
 \ov{k}&\ov{x}_k\cr 
\hline
\end{array}\,.
\end{equation}
\end{Theorem}

Here and elsewhere, for typographical convenience we have set $\ov{x}_k=x_k^{-1}$ and $\ov{y}_k=y_k^{-1}$ for $k=1,2\ldots,n$.

\section{Symplectic UASMs, GTPs and CPMs}
\label{sec-A-GT}
To make contact with other results it is necessary to introduce and relate two combinatorial constructs 
that are in bijective correspondence with unprimed symplectic shifted tableaux, namely 
certain $U$-turn alternating sign matrices (UASMs), 
and strict symplectic Gelfand-Tsetlin patterns (GTPs).


For each strict partition $\lambda$ of length $\ell(\lambda)=n$ and breadth 
$\lambda_1=m$ let ${\cal UA}^\lambda(n)$ be the set of all $2n\times m$ UASM matrices $A=(a_{ij})$,
with $i\in I$ and $j\in J=\{1,2,\ldots,m\}$,
whose matrix elements $a_{ij}$ taken from the set $\{1,0,-1\}$ are such that:\\
\begin{tabular}{rl}
UA1&the non-zero entries alternate in sign across each row and down each column;\\
UA2&the topmost non-zero entry in any column is $1$;\\
UA3&the rightmost non-zero entry in every row is $1$;\\
UA4&the sum of entries in each row and in each column is $0$ or $1$;\\
UA4&$\row_k+\row_{\ov{k}}=1$ for $k=1,2,\ldots,n$;\\
UA5&$\col_j=1$ if $j=\lambda_\ell>0$ for some $\ell$, and $=0$ otherwise, for $j=1,2,\ldots,m$,\\
\end{tabular}\\
where $\row_i$ and $\col_j$ are the sums of entries in the $i$th row and $j$th column of $A$, respectively, 
for $i\in I$ and $j\in J$.

An example of such a UASM is illustrated in the case $n=5$, $m=9$ and $\lambda=(9,7,6,2,1)$. Here
and elsewhere $\ov{1}$ has been used to denote a matrix element $-1$.
\begin{equation}
\label{eqn-UAex}
A =
\left[\!
\ASM{0.15in}{}{
{{\sc1},{\sc0},{\sc0},{\sc0},{\sc0},{\sc0},{\sc0},{\sc0},{\sc0}},
{{\sc\ov1},{\sc1},{\sc0},{\sc0},{\sc0},{\sc0},{\sc0},{\sc0},{\sc0}},
{{\sc0},{\sc0},{\sc1},{\sc0},{\sc0},{\sc0},{\sc0},{\sc0},{\sc0}},
{{\sc0},{\sc\ov1},{\sc0},{\sc1},{\sc0},{\sc0},{\sc0},{\sc0},{\sc0}},
{{\sc0},{\sc0},{\sc0},{\sc\ov1},{\sc0},{\sc1},{\sc0},{\sc0},{\sc0}},
{{\sc1},{\sc0},{\sc\ov1},{\sc1},{\sc0},{\sc0},{\sc1},{\sc0},{\sc0}},
{{\sc\ov1},{\sc1},{\sc0},{\sc\ov1},{\sc0},{\sc0},{\sc1},{\sc0},{\sc0}},
{{\sc0},{\sc0},{\sc0},{\sc0},{\sc0},{\sc0},{\sc0},{\sc1},{\sc0}},
{{\sc1},{\sc0},{\sc0},{\sc0},{\sc0},{\sc0},{\sc0},{\sc0},{\sc0}},
{{\sc0},{\sc0},{\sc0},{\sc0},{\sc0},{\sc0},{\sc0},{\sc\ov1},{\sc1}},
}
\!\right]\,.
\end{equation}

A conventional symplectic Gelfand-Tsetlin pattern $GT$ of size $2n\times n$ is an array of non-negative integers $m_{ij}$ of the form
\begin{equation}
\label{eqn-GT}
GT = \left(\!
\CPM{0.15in}{0.3in}{}{
{{m_{\ov{n}1}},{},{m_{\ov{n}2}},{},{\cdots},{},{\cdots},{},{m_{\ov{n},n}}},
{{},{m_{n1}},{},{m_{n2}},{},{\cdots},{},{m_{n,n-1}},{},{m_{nn}}},
{{},{},{\ddots},{},{\ddots},{},{},{},{\vdots}},
{{},{},{},{\ddots},{},{\ddots},{},{},{},{\vdots}},
{{},{},{},{},{m_{\ov31}},{},{m_{\ov{32}}},{},{m_{\ov33}}},
{{},{},{},{},{},{m_{31}},{},{m_{32}},{},{m_{33}}},
{{},{},{},{},{},{},{m_{\ov21}},{},{m_{\ov22}}},
{{},{},{},{},{},{},{},{m_{21}},{},{m_{22}}},
{{},{},{},{},{},{},{},{},{m_{\ov11}}},
{{},{},{},{},{},{},{},{},{},{m_{11}}},
}
\!\right)
\end{equation}
subject to the {\it betweenness conditions}
\begin{equation}\label{eqn-between}
\begin{array}{cl}
    m_{\ov{k},j} \geq m_{k,j} \geq m_{\ov{k},j+1} &\quad\hbox{for $k=1,2,\ldots,n-1$ and $j=1,2,\ldots,k$};\cr\cr
		m_{k+1,j} \geq m_{\ov{k},j} \geq m_{k+1,j+1} &\quad\hbox{for $k=1,2,\ldots,n$ and $j=1,2,\ldots,k$},
\end{array}		
\end{equation}
where $m_{\ov{k},k+1}$ is defined to be $0$ for $k=1,2,\ldots,n$.
It follows that each row is a partition.

A symplectic Gelfand-Tsetlin pattern is said to be {\it strict} if
\begin{equation}\label{eqn-strict}
    m_{kj}>m_{k,j+1} \quad\hbox{and}\quad m_{\ov{k}j}>m_{\ov{k},j+1} \quad\hbox{for $k=1,2,\ldots,n$ and $j=1,2,\ldots,k-1$},	
\end{equation}
in which case each row is a strict partition. 
For each strict partition $\lambda$ of length $\ell(\lambda)=n$ let ${\cal GT}^\lambda(n)$ be the set of all 
strict symplectic GTPs $GT$ with top row $\lambda$, that is $m_{\ov{n}j}=\lambda_j$ for $j=1,2,\ldots,n$.
We also require that $m_{kk}$ and $m_{\ov{k}k}$ are not both $0$ for any $k=1,2,\ldots,n$. This is to ensure that 
these GTPs correspond with our ST for which the first entry in the $k$th row is either $k$ or $\ov{k}$.

An example of such a strict symplectic GTP is provided in the case $n=5$ and $\lambda=(9,7,6,2,1)$ by 
\begin{equation}
\label{eqn-GTex}
GT = 
\left(\!
\CPM{0.1in}{0.1in}{}{
{{\sc9},{},{\sc7},{},{\sc6},{},{\sc2},{},{\sc1}},
{{},{\sc8},{},{\sc7},{},{\sc6},{},{\sc2},{},{\sc1}},
{{},{},{\sc8},{},{\sc7},{},{\sc6},{},{\sc2}},
{{},{},{},{\sc7},{},{\sc6},{},{\sc2},{},{\sc0}},
{{},{},{},{},{\sc6},{},{\sc4},{},{\sc1}},
{{},{},{},{},{},{\sc6},{},{\sc3},{},{\sc0}},
{{},{},{},{},{},{},{\sc4},{},{\sc3}},
{{},{},{},{},{},{},{},{\sc3},{},{\sc2}},
{{},{},{},{},{},{},{},{},{\sc2}},
{{},{},{},{},{},{},{},{},{},{\sc1}},
}
\!\right)
\end{equation}

In the case of strict symplectic Gelfand-Tsetlin patterns there are three mutually exclusive possibilities 
regarding the betweenness conditions. These may either be strict, left or right saturated, but not both.
That is to say we can define six sets of propositions:
\begin{equation}\label{eqn-BLR}
\begin{array}{lll}
  B_{kj}:=m_{kj}>m_{\ov{k-1}j}>m_{k,j+1}\,;&\quad&B_{\ov{k}j}:=m_{\ov{k}j}>m_{kj}>m_{\ov{k},j+1}\,;    \cr
	L_{kj}:=m_{kj}=m_{\ov{k-1}j}>m_{k,j+1}\,;&     &L_{\ov{k}j}:=m_{\ov{k}j}=m_{kj}>m_{\ov{k},j+1}\,;     \cr
  R_{kj}:=m_{kj}>m_{\ov{k-1}j}=m_{k,j+1}\,;&     &R_{\ov{k}j}:=m_{\ov{k}j}>m_{kj}=m_{\ov{k},j+1}\,,     \cr
\end{array}					
\end{equation}
for $1\leq j<k\leq n$ and 
\begin{equation}\label{eqn-BLRk}
\begin{array}{lll}
  B_{kk}:=m_{kk}>0\,;&\quad&B_{\ov{k}k}:=m_{\ov{k}k}>m_{kk}\,;    \cr
	L_{kk}:=m_{kk}=0\,;&     &L_{\ov{k}k}:=m_{\ov{k}k}=m_{kk}\,     \cr
	R_{kk}:=m_{kk}<0\,;&     &R_{\ov{k}k}:=m_{\ov{k}k}<m_{kk}\,     \cr	
\end{array}					
\end{equation}
for $1\leq k\leq n$,
such that for all $k=1,2,\ldots,n$ and $j=1,2,\ldots,k$
\begin{equation}\label{eqn-chi-BLR}
  \chi(B_{kj})+\chi(L_{kj})+\chi(R_{kj})=1 \quad\hbox{and}\quad \chi(B_{\ov{k}j})+\chi(L_{\ov{k}j})+\chi(R_{\ov{k}j})=1 \,,
\end{equation}
where $\chi$ is the truth function, that is to say 
$\chi(P)=1$ if $P$ is true and $=0$ if $P$ is false for any proposition $P$. 
In the case $j=k$ we have
\begin{equation}\label{eqn-chi-BLRk}
\begin{array}{c}
\chi(R_{kk})=\chi(R_{\ov{k}k})=0\,;\cr
\chi(B_{kk})+\chi(L_{kk})=1 \quad\hbox{and}\quad \chi(B_{\ov{k}k})+\chi(L_{\ov{k}k})=1 \,;\cr
\chi(L_{kk})+\chi(L_{\ov{k}k})\leq1\,,\cr 
\end{array}
\end{equation}
where the last condition corresponds to the requirement that $m_{kk}$ and $m_{\ov{k}k}$ are not both $0$.

For each strict partition $\lambda$ of length $n$ the three sets ${\cal ST}^\lambda(n)$, ${\cal UA}^\lambda(n)$ and ${\cal GT}^\lambda(n)$
are in bijective correspondence. Perhaps the simplest bijective correspondence is that between all $GT\in{\cal GT}^\lambda(n)$
and all $ST\in{\cal ST}^\lambda(n)$. For $GT=(m_{ij})$ this is defined by
\begin{equation}\label{eqn-mij-STtoGT}
  m_{ij} = \hbox{the number of entries $\leq i$ in row $j$ of $ST$} 
\end{equation}
for $i\in \{1<\ov1<2<\ov2<\cdots<n<\ov{n}\}$ and $j=1,2,\ldots,n$. 
Conversely
\begin{equation}\label{eqn-mij-GTtoST}
\begin{array}{rcl}
        \hbox{the number of entries $k$ in row $j$ of $ST$}&=& m_{kj}-m_{\ov{k-1}j}\cr
				\hbox{the number of entries $\ov{k}$ in row $j$ of $ST$}&=&m_{\ov{k}j}-m_{kj}\cr			
\end{array}
\end{equation}

To pass bijectively from $A\in{\cal UA}^\lambda(n)$ to $ST\in{\cal ST}^\lambda(n)$ one first 
constructs from $A$ a right-to-left cumulative row sum matrix, $\row(A)$. The non-zero entries $1$ in column $j$ of $\row(A)$ then specify 
by their row number an entry $k$ or $\ov{k}$ on the $j$th diagonal of $ST$.
To complete the map from $\row(A)$ to $ST$, these row numbers are arranged in strictly increasing order down each diagonal 
from top-left to bottom-right.

Similarly, to pass bijectively from $A\in{\cal UA}^\lambda(n)$ to $GT\in{\cal GT}^\lambda(n)$ one first constructs from 
$A$ a top-to-bottom cumulative column sum matrix, $\col(A)$. The non-zero entries $1$ in row $i$ of $\col(A)$,
counted from top to bottom then specify 
by their column number an entry $j$ in row $i$ of $GT$
counted from bottom to top. These column numbers are arranged 
in strictly increasing order across each row from left to right starting on the main diagonal with spaces
allowing entries in successive rows to lie between one another. 

That each of these simple, invertible maps is a bijection is not immediately obvious, but is not hard to verify.  

To illustrate the above bijections with our running example, taken from (\ref{eqn-ST}), the maps between $A$, $\row(A)$ and $ST$ 
are shown below across the top row of (\ref{eqn-ASCG}), and those between $A$, $\col(A)$ and $GT$ are shown below
in the left hand column of (\ref{eqn-ASCG}). 

\begin{center}
\begin{equation}\label{eqn-ASCG}
\begin{array}{l}
A =
\left[\!
\ASM{0.15in}{}{
{{\sc1},{\sc0},{\sc0},{\sc0},{\sc0},{\sc0},{\sc0},{\sc0},{\sc0}},
{{\sc\ov1},{\sc1},{\sc0},{\sc0},{\sc0},{\sc0},{\sc0},{\sc0},{\sc0}},
{{\sc0},{\sc0},{\sc1},{\sc0},{\sc0},{\sc0},{\sc0},{\sc0},{\sc0}},
{{\sc0},{\sc\ov1},{\sc0},{\sc1},{\sc0},{\sc0},{\sc0},{\sc0},{\sc0}},
{{\sc0},{\sc0},{\sc0},{\sc\ov1},{\sc0},{\sc1},{\sc0},{\sc0},{\sc0}},
{{\sc1},{\sc0},{\sc\ov1},{\sc1},{\sc0},{\sc0},{\sc0},{\sc0},{\sc0}},
{{\sc\ov1},{\sc1},{\sc0},{\sc\ov1},{\sc0},{\sc0},{\sc1},{\sc0},{\sc0}},
{{\sc0},{\sc0},{\sc0},{\sc0},{\sc0},{\sc0},{\sc0},{\sc1},{\sc0}},
{{\sc1},{\sc0},{\sc0},{\sc0},{\sc0},{\sc0},{\sc0},{\sc0},{\sc0}},
{{\sc0},{\sc0},{\sc0},{\sc0},{\sc0},{\sc0},{\sc0},{\sc\ov1},{\sc1}},
}
\!\right]
\CPM{0.15in}{0.1in}{}{
{{\sc1}},
{{\sc\ov1}},
{{\sc2}},
{{\sc\ov2}},
{{\sc3}},
{{\sc\ov3}},
{{\sc4}},
{{\sc\ov4}},
{{\sc5}},
{{\sc\ov5}},
}
\Leftrightarrow
\row(A) =
\left[\!
\ASM{0.15in}{}{
{{\sc1},{\sc0},{\sc0},{\sc0},{\sc0},{\sc0},{\sc0},{\sc0},{\sc0}},
{{\sc0},{\sc1},{\sc0},{\sc0},{\sc0},{\sc0},{\sc0},{\sc0},{\sc0}},
{{\sc1},{\sc1},{\sc1},{\sc0},{\sc0},{\sc0},{\sc0},{\sc0},{\sc0}},
{{\sc0},{\sc0},{\sc1},{\sc1},{\sc0},{\sc0},{\sc0},{\sc0},{\sc0}},
{{\sc0},{\sc0},{\sc0},{\sc0},{\sc1},{\sc1},{\sc0},{\sc0},{\sc0}},
{{\sc1},{\sc0},{\sc0},{\sc1},{\sc0},{\sc0},{\sc0},{\sc0},{\sc0}},
{{\sc0},{\sc1},{\sc0},{\sc0},{\sc1},{\sc1},{\sc1},{\sc0},{\sc0}},
{{\sc1},{\sc1},{\sc1},{\sc1},{\sc1},{\sc1},{\sc1},{\sc1},{\sc0}},
{{\sc1},{\sc0},{\sc0},{\sc0},{\sc0},{\sc0},{\sc0},{\sc0},{\sc0}},
{{\sc0},{\sc0},{\sc0},{\sc0},{\sc0},{\sc0},{\sc0},{\sc0},{\sc1}},
}
\!\right]

\Leftrightarrow
ST =\!\!\!\!
\SYT{0.15in}{}{
 {{\sc1},{\sc\ov1},{\sc2},{\sc\ov2},{\sc3},{\sc3},{\sc4},{\sc\ov4},{\sc\ov5}},
 {{\sc2},{\sc2},{\sc\ov2},{\sc\ov3},{\sc4},{\sc4},{\sc\ov4}},
 {{\sc\ov3},{\sc4},{\sc\ov4},{\sc\ov4},{\sc\ov4},{\sc\ov4}},
 {{\sc\ov4},{\sc\ov4}},
 {{\sc5}}
}
\cr\cr
{\hskip7em}\Updownarrow{\hskip7em}\searrow{\hskip-1em}\nwarrow
\cr\cr
\col(A) =
\left[\!
\ASM{0.15in}{}{
{{\sc1},{\sc0},{\sc0},{\sc0},{\sc0},{\sc0},{\sc0},{\sc0},{\sc0}},
{{\sc0},{\sc1},{\sc0},{\sc0},{\sc0},{\sc0},{\sc0},{\sc0},{\sc0}},
{{\sc0},{\sc1},{\sc1},{\sc0},{\sc0},{\sc0},{\sc0},{\sc0},{\sc0}},
{{\sc0},{\sc0},{\sc1},{\sc1},{\sc0},{\sc0},{\sc0},{\sc0},{\sc0}},
{{\sc0},{\sc0},{\sc1},{\sc0},{\sc0},{\sc1},{\sc0},{\sc0},{\sc0}},
{{\sc1},{\sc0},{\sc0},{\sc1},{\sc0},{\sc1},{\sc0},{\sc0},{\sc0}},
{{\sc0},{\sc1},{\sc0},{\sc0},{\sc0},{\sc1},{\sc1},{\sc0},{\sc0}},
{{\sc0},{\sc1},{\sc0},{\sc0},{\sc0},{\sc1},{\sc1},{\sc1},{\sc0}},
{{\sc1},{\sc1},{\sc0},{\sc0},{\sc0},{\sc1},{\sc1},{\sc1},{\sc0}},
{{\sc1},{\sc1},{\sc0},{\sc0},{\sc0},{\sc1},{\sc1},{\sc0},{\sc1}},
}
\!\right]
~~~
C(A) =
\left[\!
\CPM{0.15in}{0.3in}{}{
{{\sc{\WE}},{\sc{\SE}},{\sc{\SE}},{\sc{\SE}},{\sc{\SE}},{\sc{\SE}},{\sc{\SE}},{\sc{\SE}},{\sc{\SE}}},
{{\sc{\NS}},{\sc{\WE}},{\sc{\SE}},{\sc{\SE}},{\sc{\SE}},{\sc{\SE}},{\sc{\SE}},{\sc{\SE}},{\sc{\SE}}},
{{\sc{\SW}},{\sc{\NW}},{\sc{\WE}},{\sc{\SE}},{\sc{\SE}},{\sc{\SE}},{\sc{\SE}},{\sc{\SE}},{\sc{\SE}}},
{{\sc{\SE}},{\sc{\NS}},{\sc{\NW}},{\sc{\WE}},{\sc{\SE}},{\sc{\SE}},{\sc{\SE}},{\sc{\SE}},{\sc{\SE}}},
{{\sc{\SE}},{\sc{\SE}},{\sc{\NE}},{\sc{\NS}},{\sc{\SW}},{\sc{\WE}},{\sc{\SE}},{\sc{\SE}},{\sc{\SE}}},
{{\sc{\WE}},{\sc{\SE}},{\sc{\NS}},{\sc{\WE}},{\sc{\SE}},{\sc{\NE}},{\sc{\SE}},{\sc{\SE}},{\sc{\SE}}},
{{\sc{\NS}},{\sc{\WE}},{\sc{\SE}},{\sc{\NS}},{\sc{\SW}},{\sc{\NW}},{\sc{\WE}},{\sc{\SE}},{\sc{\SE}}},
{{\sc{\SW}},{\sc{\NW}},{\sc{\SW}},{\sc{\SW}},{\sc{\SW}},{\sc{\NW}},{\sc{\NW}},{\sc{\WE}},{\sc{\SE}}},
{{\sc{\WE}},{\sc{\NE}},{\sc{\SE}},{\sc{\SE}},{\sc{\SE}},{\sc{\NE}},{\sc{\NE}},{\sc{\NE}},{\sc{\SE}}},
{{\sc{\NE}},{\sc{\NE}},{\sc{\SE}},{\sc{\SE}},{\sc{\SE}},{\sc{\NE}},{\sc{\NE}},{\sc{\NS}},{\sc{\WE}}},
}
\!\right]
\CPM{0.15in}{0.1in}{}{
{{\sc1}},
{{\sc\ov1}},
{{\sc2}},
{{\sc\ov2}},
{{\sc3}},
{{\sc\ov3}},
{{\sc4}},
{{\sc\ov4}},
{{\sc5}},
{{\sc\ov5}},
}
\cr\cr
{\hskip7em}\Updownarrow
\cr\cr
GT = 
\left(\!
\CPM{0.15in}{0.1in}{}{
{{\sc9},{},{\sc7},{},{\sc6},{},{\sc2},{},{\sc1}},
{{},{\sc8},{},{\sc7},{},{\sc6},{},{\sc2},{},{\sc1}},
{{},{},{\sc8},{},{\sc7},{},{\sc6},{},{\sc2}},
{{},{},{},{\sc7},{},{\sc6},{},{\sc2},{},{\sc0}},
{{},{},{},{},{\sc6},{},{\sc4},{},{\sc1}},
{{},{},{},{},{},{\sc6},{},{\sc3},{},{\sc0}},
{{},{},{},{},{},{},{\sc4},{},{\sc3}},
{{},{},{},{},{},{},{},{\sc3},{},{\sc2}},
{{},{},{},{},{},{},{},{},{\sc2}},
{{},{},{},{},{},{},{},{},{},{\sc1}},
}
\!\right)
\CPM{0.15in}{0.1in}{}{
{{\sc\ov5}},
{{\sc5}},
{{\sc\ov4}},
{{\sc4}},
{{\sc\ov3}},
{{\sc3}},
{{\sc\ov2}},
{{\sc2}},
{{\sc\ov1}},
{{\sc1}},
}
\cr
\end{array}
\end{equation}
\end{center}

In this diagram we have also illustrated a $U$-turn compass point matrix (CPM) $C(A)$ that is obtained from $A$
by mapping the entries $1$ and $\ov1$ in $A$ to $\WE$ and $\NS$, respectively, 
and the entries $0$ in $A$ to one or other of $\NE$, $\SE$, $\NW$ or $\SW$ 
in accordance with the arrangements of the nearest non-zero  
neighbours of the $0$ in the four compass point directions, as specified in the following tabulation. 
In the absence of any non-zero element of $A$ in the directions $N$ or $E$
such a missing element is deemed to be $\ov1$.
\begin{equation}\label{eqn-ASM-CPM}
\begin{array}{|c|c|c|c|c|c|c|}
\hline
\hbox{UASM entry $a_{ij}^{\phantom{1}}$} & 1 &\ov1 & 0 & 0 &0 &0 \cr
\hline
&&&&&&\cr
                &\begin{array}{c}\ov1\cr \ov1~~{\bf1}~~\ov1\cr\ov1\end{array}
								&\begin{array}{c} 1\cr 1~~{\bf\ov1}~~1\cr 1\end{array}
								&\begin{array}{c} 1\cr 1~~{\bf0}~~\ov1\cr\ov1\end{array}
								&\begin{array}{c}\ov1\cr 1~~{\bf0}~~\ov1\cr 1\end{array}
								&\begin{array}{c} 1\cr \ov1~~{\bf0}~~1\cr\ov1\end{array}
								&\begin{array}{c}\ov1\cr \ov1~~{\bf0}~~1\cr 1\end{array}\cr
&&&&&&\cr								
\hline
\hbox{CPM entry $c_{ij}^{\phantom{1}}$} & \WE &\NS & \NE & \SE & \NW & \SW \cr
\hline
\end{array}
\end{equation}
where once again we have adopted a notation whereby $\ov1$ signifies an entry $-1$, and the notation of~\cite{HK-bij} 
has been altered by interchanging both $N$ with $S$, and $W$ with $E$. 

We define the set ${\cal UC}^\lambda(n)$ of CPMs corresponding to $A\in{\cal UA}^\lambda(n)$ 
to be those matrices $C(A)$ obtained in the above manner. They have an important role to play in establishing the equivalence 
of various possible ways of weighting entries in symplectic shifted tableaux and in strict symplectic Gelfand-Tsetlin patterns.

\section{Alternative forms of the symplectic Tokuyama identity}
\label{sec-Tok-A-GT}
In order to establish corollaries of Theorem~\ref{thm-main} within the context of the above combinatorial objects
it is merely necessary to replace the sum over $ST\in{\cal ST}^\lambda(\x/\y)$ by sums over $A\in{\cal UA}^\lambda$
and $GT\in{\cal GT}^\lambda$, along with appropriate identifications of $\wgt(c_{ij})$ and $\wgt(m_{ij})$. 

The simplest case is that of ${\cal UA}^\lambda(n)$. The non-zero entries in $ST$ arising from the $1$s appearing in $\row(A)$
are associated with the entries $\WE$, $\SW$ and $\NW$ in $C(A)$. These entries in $ST$ are $k$ or $\ov{k}$ according as the entry 
in $C(A)$ lies in its row $k$ or $\ov{k}$, respectively. The fact that the $k$th entry on the main diagonal of $ST$
is either $k$ or $\ov{k}$ is a consequence of the $U$-turn nature of $A$ that ensures that either $c_{k1}$ or $c_{\ov{k}1}$  belongs to 
the set $\{\WE,\SW,\NW\}$, but not both. Moreover the entries $\SW$ and $\NW$ in $C(A)$ are associated with
identical pairs of entries in $ST$ that are horizontal and vertical, respectively. It follows from the 
tabulation of weights given in (\ref{eqn-swgt-twgt}) for all entries $s_{ij}$ in $ST$ that 
Theorem~\ref{thm-main} may be rewritten in the form:
\begin{Corollary}\label{cor-A}
Let $\mu$ be a partition of length $\ell(\mu)\leq n$ and $\delta= (n, n-1, \ldots, 1)$, so that
$\lambda=\mu+\delta$ is a strict partition of length $\ell(\lambda)=n$. 
Then for all $\x=(x_1,x_2,\ldots,x_n)$ and $\y=(y_1,y_2,\ldots,y_n)$
\begin{equation}\label{eqn-Tok-A}
     \sum_{A\in{\cal UA}^\lambda(n)} \wgt(A)  = \prod_{1\leq i\leq j \leq n} (x_i + y_j)(1+x_i^{-1}y_j^{-1})\ \ sp_{\mu} (\x)\,,
\end{equation}
where 
\begin{equation}\label{eqn-wgtAC}
      \wgt(A) = \prod_{i\in I}\prod_{j=1}^m \wgt(c_{ij})\,.
\end{equation}
with
\begin{equation}\label{eqn-wgtC}
\begin{array}{|c|c|c|c|c|c|c|}
\hline
 c_{ij}~~1\leq j\leq m&\WE&\NS&\NE&\SE&\NW&\SW\cr
\hline
i=k^{\phantom1}&x_k+y_k&1&1&1&y_k&x_k\cr
\hline
i=\ov{k}^{\phantom1}&\ov{x}_k+\ov{y}_k&1&1&1&\ov{y}_k&\ov{x}_k\cr
\hline
\end{array}
\end{equation}
\end{Corollary}

This serves to correct the corresponding result given in Corollary~5.4 of~\cite{HK-bij}. The correction
consists of changing $i<j$ to $i\leq j$ in the first line of equation (5.94),
and of replacing $\NS_k(A)$ and $\NS_{\ov{k}}(A)$ to $\WE_k(A)$ and $\WE_{\ov{k}}(A)$, respectively, in the third line of equation (5.94). 
This is necessary to obtain the factorisation as claimed.

It might be noted here that the alternating sign nature of $A$ implies that 
if $\#XY_{i}$ is the number of entries $XY$ in row $i$ of $C(A)$ for all pairs of compass point
directions $X$ and $Y$, then
\begin{equation}\label{eqn-WE-NS}
   \#\WE_i=\#\NS_i+\chi(P_i) \quad\hbox{with}\quad P_i:= c_{i1}\in\{\WE,\SW,\NW\}\,.
\end{equation}
With this notation the $U$-turn nature of $A$ is reflected in the fact that 
\begin{equation}
\chi(P_{k})+\chi(P_{\ov{k}})=1\,.
\end{equation}

Corollary~\ref{cor-A} then still applies with the weighting specified by:
\begin{equation}\label{eqn-wgtC-alt}
\begin{array}{|c|c|c|c|c|c|c|}
\hline
 c_{ij}~~1\leq j\leq m&WE&NS&NE&SE&NW&SW\cr
\hline
i=k,j=1^{\phantom{1}}&x_k+y_k&1&1&1&x_k+y_k&x_k+y_k\cr
\hline
i=\ov{k},j=1^{\phantom{k^2}}&\ov{x}_k+\ov{y}_k&1&1&1&\ov{x}_k+\ov{y}_k&\ov{x}_k+\ov{y}_k\cr
\hline
i=k,j\geq 1^{\phantom{1}}&1&x_k+y_k&1&1&y_k&x_k\cr
\hline
i=\ov{k},j\geq 1^{\phantom{k^2}}&1&\ov{x}_k+\ov{y}_k&1&1&\ov{y}_k&\ov{x}_k\cr
\hline
\end{array}
\end{equation}
where associating factors $x_k+y_k$ and $\ov{x}_k+\ov{y}_k$ with all the entries $\NS$ rather than
$\WE$ in $C(A)$ is compensated for by the inclusion of the additional $j=1$ weights dictated by (\ref{eqn-WE-NS}).

\medskip

Turning to Gelfand-Tsetlin patterns, 
the analog of Corollary~\ref{cor-A} takes the form:

\begin{Corollary}\label{cor-TokGT}
Let $\mu$ be a partition of length $\ell(\mu)\leq n$ and $\delta= (n, n-1, \ldots, 1)$, so that
$\lambda=\mu+\delta$ is a strict partition of length $\ell(\lambda)=n$. 
Then for all $\x=(x_1,x_2,\ldots,x_n)$ and $\y=(y_1,y_2,\ldots,y_n)$
\begin{equation}\label{eqn-Tok-AGT}
     \sum_{GT\in{\cal GT}^\lambda(n)} \wgt(GT)  = \prod_{1\leq i\leq j \leq n} (x_i + y_j)(1+x_i^{-1}y_j^{-1})\ \ sp_{\mu} (\x)\,,
\end{equation}
where 
\begin{equation}\label{eqn-wgtGT}
\begin{array}{rcl}
  \wgt(GT)
	        &=& \ds\prod_{k=1}^n\ \prod_{j=1}^{k}\ 
	                  \left(\,\chi(B_{kj}) (x_k+y_k) + \chi(L_{kj})x_k + \chi(R_{kj})y_k\,\right)\ x_k^{m_{kj}-m_{\ov{k-1}j}-1}\cr
					 &\times&\ds\prod_{k=1}^n\ \prod_{j=1}^{k}\ 
					         \left(\,\chi(B_{\ov{k}j}) (\ov{x}_k+\ov{y}_k) + \chi(L_{\ov{k}j})\ov{x}_k+ \chi(R_{\ov{k}j})\ov{y}_k \,\right)\ \ov{x}_k^{m_{\ov{k}j}-m_{kj}-1} \,,\cr
\end{array}					
\end{equation}
with $m_{\ov{k-1}k}=0$ for $k=1,2,\ldots,n$.
\end{Corollary}

\noindent{\bf Proof}:
Let $ST$ and $A$ be the elements of ${\cal ST}^\lambda(n)$ and ${\cal UA}^\lambda(n)$ in bijective correspondence 
with a given $GT\in{\cal GT}^\lambda(n)$, and let $C(A)$ be the corresponding compass point matrix. 
The relationship (\ref{eqn-mij-GTtoST}) between $ST$ and $GT$ implies that the sequence of entries $k$ in row $j$ of $ST$
has length $\ell=m_{kj}-m_{\ov{k-1}j}$. In the case $j<k$ if $\ell>0$ the leftmost of these entries will be either 
the beginning of a continuous strip of $k$s in row $j$, or the continuation of a continuous strip 
of $k$'s from row $j+1$ with an identical entry $k$ immediately beneath it. All other entries $k$ in row
$j$ will have an identical entry $k$ immediately to its left. These two cases correspond
to the existence of an entry $\NS$ or $\NW$ in row $k$ of $C(A)$ together with a sequence of $\ell-1$ entries $\SW$.
However, the map from $A$ to $GT$ by way of $\col(A)$ is such that in the
notation of (\ref{eqn-BLR}) the betweenness conditions $\chi(B_{kj})=1$, $\chi(R_{kj})=1$ 
and $\chi(L_{kj})=1$ correspond to the existence of an entry $\NS$, $\NW$ and $\NE$, respectively, 
in row $k$ of $C(A)$. It follows from the weighting (\ref{eqn-wgtC-alt}) of $C(A)$ that the weights
attached to $\chi(B_{kj})$ and $\chi(R_{kj})$ must be $(x_k+y_k)x_k^{\ell-1}$ and $y_k\,x_k^{\ell-1}$, 
respectively, as indicated in the top-line of (\ref{eqn-wgtGT}). 
On the otherhand, if $\ell=0$ there are no entries $k$ in row $j$ so the contribution to the
weight should simply be a factor $1$. This is reflected in the fact that in this case $\chi(L_{kj})=1$ and 
the weight attached to it in (\ref{eqn-wgtGT}) is indeed just $x_k\,x_k^{0-1}=1$, as required.

The case $j=k$ is slightly different. The sequence of entries $k$ in row $k$ of $ST$
now has length $\ell=m_{kk}$, which is of course equal to $m_{kk}-m_{\ov{k-1},k}$ since $m_{\ov{k-1},k}$
has been defined to be $0$ by hypothesis. For $\ell>0$ the leftmost entry $k$ necessarily lies on the
main diagonal of $ST$ and carries a weight $x_k+y_k$ with the remaining $\ell-1$ entries $k$ in row $k$ all carrying
a weight $x_k$. This is correctly codified in (\ref{eqn-wgtGT}) since in this case we have $m_{kk}>0$
so that $\chi(B_{kk})=1$ and $\chi(L_{kk})=\chi(R_{kk})=0$. If $\ell=0$ there is no contribution to 
$\wgt(GT)$ and this is again in accord with (\ref{eqn-wgtGT}) since $\ell=m_{kk}=0$ implies that
$\chi(L_{kk})=1$ and $\chi(B_{kk})=\chi(R_{kk})=0$, with the weight attached to $\chi(L_{kk})$ being
$x_k\,x_k^{0-1}=1$, as required.

Sequences of entries $\ov{k}$ of length $\ell=m_{\ov{k}k}-m_{kk}$ in row $j$ of $ST$ can be dealt 
with in a similar manner so as to confirm the validity of the second line of (\ref{eqn-wgtGT}). 
If $\ell>0$ and $j<k$ then the argument is exactly the same as before with $k$ replaced by $\ov{k}$. 
If $\ell>0$ and $j=k$ then the sequence of entries $\ov{k}$ in row $k$ of $ST$ only reach the main diagonal if $m_{kk}=0$,
but, whether or not this is the case, they carry a total weight of $(\ov{x}_k+\ov{y}_k)x_k^{\ell-1}$. This is 
once again in accord with (\ref{eqn-wgtGT}) since $\chi(B_{\ov{k}k})=1$ and $\chi(L_{\ov{k}k})=\chi(R_{\ov{k}k})=0$.
If $\ell=0$ we have both $\chi(B_{\ov{k}j})=0$ and $\chi(R_{\ov{k}j})=0$ for all $j$, leaving $\chi(L_{\ov{k}j})=1$.
However, this leads as usual to a weight contribution of $\ov{x}_k\,\ov{x}_k^{0-1}=1$, again just as required.
\qed

\medskip

It might be worth noting that the weight of the various combinatorial entities displayed in (\ref{eqn-ASCG}) are given by
\begin{equation}\label{eqn-example}
  \wgt(ST)=\wgt(A)=\wgt(GT) = \frac{(x_1+y_1)^2(x_2+y_2)^2(x_3+y_3)^3(x_4+y_4)^3(x_5+y_5)^2}{x_1x_3x_4^4x_5\,y_1y_2y_3^2y_4^3y_5}
\end{equation}

The tableau from (\ref{eqn-ASCG}) can also be written so the weights associated with each box are clearly visible:

\begin{equation}
\wideSYT{0.2in}{0.6in}{}{
 {{x_1+y_1},{\ov{x}_1+ \ov{y}_1},{y_2},{\ov{y}_2},{x_3+ y_3},{x_3},{y_4},{\ov{y}_4},{\ov{x}_5+ \ov{y}_5}},
 {{x_2+y_2},{x_2},{\ov{x}_2+ \ov{y}_2},{\ov{x}_3+ \ov{y}_3},{x_4+ y_4},{x_4},{\ov{y}_4}},
 {{\ov{x}_3+ \ov{y}_3},{x_4+y_4},{\ov{y}_4},{\ov{x}_4},{\ov{x}_4},{\ov{x}_4}},
 {{\ov{x}_4+ \ov{y}_4},{\ov{x}_4}},
 {{x_5+ y_5}}
}
\end{equation}

\section{The specialisation $y_k=qx_k$ for all $k$}
\label{sec-yqx}
The specialisation $y_k=qx_k$ for $k=1,2\ldots,n$ for any non-zero parameter $q$ leads to a considerable simplification 
of all our identities.

\begin{Corollary}\label{cor-main-q}
Let $\mu$ be a partition of length $\ell(\mu)\leq n$ and $\delta= (n, n-1, \ldots, 1)$, so that
$\lambda=\mu+\delta$ is a strict partition of length $\ell(\lambda)=n$. 
Then for all $\x=(x_1,x_2,\ldots,x_n)$ and $q$
\begin{equation}\label{eqn-Tok-ST-q}
    \sum_{ST\in{\cal ST}^\lambda(\x)}\ \wgt(ST) = \prod_{1\leq i\leq j \leq n} (x_i + qx_j)(1+q^{-1}x_i^{-1}x_j^{-1})\ \ sp_\mu(\x)\,,
\end{equation}
where
\begin{equation}
   \wgt(ST)=\prod_{(i,j)\in{SF^\lambda}}\!\!\!\wgt(s_{ij}) 
\end{equation}
with 
\begin{equation}\label{eqn-swgt-q}
\begin{array}{||c|c|l||c|c|l||}
\hline
s_{ij}&\wgt(s_{ij})&&s_{ij}&\wgt(s_{ij})&\cr
\hline
 k&x_k &\hbox{if $s_{i,j-1}=k$}&\ov{k}^{\phantom{1}}&\ov{x}_k &\hbox{if  $s_{i,j-1}=\ov{k}$}\cr
  &q\,x_k &\hbox{if $s_{i-1,j}=k$} & &(1/q)\ov{x}_k &\hbox{if $s_{i-1,j}=\ov{k}$}\cr
	&(1+q)x_k&\hbox{otherwise} & &(1+1/q)\ov{x}_k &\hbox{otherwise}\cr
\hline
\end{array}
\end{equation}
\end{Corollary}

Before proceeding to the next corollary it is convenient to introduce a small lemma
\begin{Lemma}
\label{lem-NE-NW}
Let $C(A)$ be the compass point matrix corresponding to the UASM $A\in{\cal A}^\lambda$, and
Let $\#\XY_i$ be the number of entries $\XY$ in the $i$th row of $C(A)$ for $i=k$ or $\ov{k}$
and $k=1,2,\ldots,n$. Then
\begin{equation}\label{eqn-lem}
\begin{array}{rcl}
\#\NS_{k}+\#\NW_{k}+\#\NE_{k} &=&k-1\,,\cr
\#\WE_{\ov{k}}+\#\NW_{\ov{k}}+\#\NE_{\ov{k}}&=&k\,.\cr
\end{array}
\end{equation}  
\end{Lemma}

\noindent{\bf Proof}:~~
The first result follows from
\begin{equation}
\#\NS_{k}+\#\NW_{k}+\#\NE_{k} = \sum_{i=1}^{\ov{k-1}} \sum_{j=1}^m a_{ij} = k-1\,,
\end{equation}
where the first step follows from the fact that the tabulation of (\ref{eqn-ASM-CPM}) implies
that the column sum in $A$ above the position of each entry $\XY$ in row $k$ of $C$ 
is $1$ or $0$ according as $\XY$ is or is not in $\{\NS,\NW,\NE\}$. The second step is a consequence 
of the fact that the row sums $\row_\ell$ and $\row_{\ov{\ell}}$ for consecutive pairs of rows $i=\ell$ and $i=\ov{\ell}$
are such that $\row_\ell+\row_{\ov{\ell}}=1$ for $\ell=1,2,\ldots,k-1$.

In the second case we have
\begin{equation}
\#\NS_{\ov{k}}+\#\NW_{\ov{k}}+\#\NE_{\ov{k}} 
=\sum_{i=1}^{k} \sum_{j=1}^m a_{ij} = k-1+\chi(P_k) = k-\chi(P_{\ov{k}})\,,
\end{equation}
where the first step is essentially the same as before and the second relies on the fact
that the row sum $\row_k$ is $1$ or $0$ according as $c_{k1}$ is or is not $\in\{\WE,\SW,\NW\}$,
that is to say according as $\chi(P_k)$, as defined in (\ref{eqn-WE-NS}), is $1$ or $0$.
The third step just uses the fact that $\chi(P_k)+\chi(P_{\ov{k}})=1$. Finally,
if one combines this result with the identity (\ref{eqn-WE-NS}) one obtains (\ref{eqn-lem}),
as required.
\qed
\medskip

Now we are in a position to state the following
\begin{Corollary}\label{cor-TokAq}
Let $\mu$ be a partition of length $\ell(\mu)\leq n$ and $\delta= (n, n-1, \ldots, 1)$, so that
$\lambda=\mu+\delta$ is a strict partition of length $\ell(\lambda)=n$. 
Then for all $\x=(x_1,x_2,\ldots,x_n)$ and $q$
\begin{equation}\label{eqn-Tok-Aq}
     \sum_{A\in{\cal UA}^\lambda(n)} \wgt(A)  = \prod_{1\leq i\leq j\leq n}\!\!\!(x_i+q\,x_j)(1+q^{-1}x_i^{-1}x_j^{-1})\ \ sp_\mu(\x)\,,
\end{equation}
where 
\begin{equation}\label{eqn-wgtACq}
      \wgt(A) =  c_0\ \prod_{i=1}^{2n}\prod_{j=1}^m \wgt(c_{ij})\,,
\end{equation}
with $c_0=1$ and
\begin{equation}\label{eqn-wgtCqa}
\begin{array}{|c|c|c|c|c|c|c|}
\hline
c_{ij}~~1\leq j\leq m&WE&NS&NE&SE&NW&SW\cr
\hline
i=2k-1^{\phantom1}&(1+q)x_k&1&1&1&qx_k&x_k\cr
\hline
i=2k^{\phantom1}&(1+1/q)\ov{x}_k&1&1&1&(1/q)\ov{x}_k&\ov{x}_k\cr
\hline
\end{array}
\end{equation}
or equivalently $c_0=(1+q)/q^{n(n+1)/2}$ and
\begin{equation}\label{eqn-wgtCqb}
\begin{array}{|c|c|c|c|c|c|c|}
\hline
c_{ij}~~1\leq j\leq m&WE&NS&NE&SE&NW&SW\cr
\hline
i=2k-1^{\phantom1}&x_k&(1+q)&1&1&qx_k&x_k\cr
\hline
i=2k^{\phantom1}&\ov{x}_k&(1+q)&q&1&\ov{x}_k&\ov{x}_k\cr
\hline
\end{array}
\end{equation}
\end{Corollary}

\noindent{\bf Proof}:~~
The first form of the weighting with the overall factor $c_0=1$ is a direct consequence of 
Corollary~\ref{cor-A} with $y_k$ set equal to $q\,x_k$ for all $k$. The second form arises
by writing each factor $(1+1/q)$ in the $\WE$ column of the $\ov{k}$ row in the form $(1/q)(1+q)$,
and the factor $1$ in the $\NE$ column of the $\ov{k}$ row in the form $(1/q)q$. This allows
one to extract a factor $(1/q)^k$ as a result of Lemma~\ref{lem-NE-NW}. Taking the product over $k$ from
$1$ to $n$ gives an overall factor of $1/q^{n(n+1)/2}$ and leaves contributions $(1+q)$ and $q$ in 
the $\WE$ and $\NE$ columns. Finally, one can use the identity (\ref{eqn-WE-NS}) to move the
factors $(1+q)$ from the $\WE$ to the $\NS$ column in rows $k$ and $\ov{k}$, leading
to an overall factor of $(1-q)$ to the power $\sum_{k=1}^n(\chi(P_k)+\chi(P_{\ov{k}}))=n$.
\qed.
\medskip

Similarly, in the terms of symplectic Gelfand-Tsetlin patterns we have:
 
\begin{Corollary}\label{cor-TokGTq}
Let $\mu$ be a partition of length $\ell(\mu)\leq n$ and $\delta= (n, n-1, \ldots, 1)$, so that
$\lambda=\mu+\delta$ is a strict partition of length $\ell(\lambda)=n$. 
Then for all $\x=(x_1,x_2,\ldots,x_n)$ and $q$
\begin{equation}\label{eqn-Tok-GTq}
     \sum_{GT\in{\cal GT}^\lambda(n)} \wgt(GT)  = \prod_{1\leq i\leq j \leq n} (x_i + q\,x_j)(1+q^{-1}x_i^{-1}x_j^{-1})\ \ sp_{\mu} (\x)\,,
\end{equation}
where 
\begin{equation}\label{eqn-wgtGTq}
 \wgt(GT) = \prod_{k=1}^n\ \prod_{j=1}^{k}\ (1+q)^{\chi(B_{kj})+\chi(B_{\ov{k}j})}\ q^{\chi(R_{kj})+\chi(L_{\ov{k}j})-1}\  \x^{\xwgt(GT)}
\end{equation}
with 
\begin{equation}
      \x^{\xwgt(GT)} = \prod_{k=1}^n\ \prod_{j=1}^k x_k^{ 2m_{kj}-m_{\ov{k}j}-m_{\ov{k-1}j}} \,,
\end{equation}
and $m_{\ov{k-1}k}=0$ for $k=1,2,\ldots,n$.
\end{Corollary}

\noindent{\bf Proof}:
Setting $y_k=q\,x_k$ for all $k$ in (\ref{eqn-wgtGT}) gives
\begin{equation}\label{eqn-wgtGTqa}
\begin{array}{rcl}
  \wgt(GT)&=&\ds\prod_{k=1}^n\ \prod_{j=1}^{k}\ 
	                  \left(\,\chi(B_{kj}) (1+q)x_k + \chi(L_{kj})x_k + \chi(R_{kj})q\,x_k\,\right)\ x_k^{m_{kj}-m_{\ov{k-1}j}-1}\cr
					&\times&\ds\prod_{k=1}^n\ \prod_{j=1}^{k}\ 
					         \left(\,\chi(B_{\ov{k}j}) (1+q)\ov{x}_k + \chi(L_{\ov{k}j})q\ov{x}_k+ \chi(R_{\ov{k}j})\ov{x}_k \,\right)\,q^{-1}\,\ov{x}_k^{m_{\ov{k}j}-m_{kj}-1} \cr
\end{array}					
\end{equation}
From this one can see immediately that the $\x$ dependence is given by
\begin{equation}\label{eqn-x-dependence}
               \prod_{k=1}^n\ \prod_{j=1}^k  x_k^{ 2m_{kj}-m_{\ov{k}j}-m_{\ov{k-1}j}}\,.
\end{equation}
and the $q$ dependence by 
\begin{equation}\label{eqn-q-dependence}
              \prod_{k=1}^n\ \prod_{j=1}^k  (1+q)^{\chi(B_{kj})+\chi(B_{\ov{k}j})}  q^{\chi(R_{kj})+\chi(L_{\ov{k}j})-1} \,,
\end{equation}
as required.
\qed
\medskip

This result may be expressed in the alternative form first obtained by Gupta {\em et al.}~\cite{GRVP}:
 
\begin{Corollary}\label{cor-TokGTq-final}
Let $\mu$ be a partition of length $\ell(\mu)\leq n$ and $\delta= (n, n-1, \ldots, 1)$, so that
$\lambda=\mu+\delta$ is a strict partition of length $\ell(\lambda)=n$. 
Then for all $\x=(x_1,x_2,\ldots,x_n)$ and $q$
\begin{equation}\label{eqn-Tok-GTqq}
     \sum_{GT\in{\cal GT}^\lambda(n)} \qxwgt(GT)  = \prod_{i=1}^n (q\,x_i + x_i^{-1}) \ \prod_{1\leq i<j \leq n} (x_i + q\,x_j)(q+x_i^{-1}x_j^{-1})\ \ sp_{\mu} (\x)\,,
\end{equation}
where 
\begin{equation}\label{eqn-wgtGTq-final}
   \qxwgt(GT) = (1+q)^{B(GT)} \ q^{R_o(GT)+L_e(GT)}\ \x^{\xwgt(GT)}			
\end{equation}
with
\begin{equation}\label{eqn-BRoLe}
\begin{array}{rcl}
     B(GT)&=& \#\{(k,j)\,|\, m_{k+1,j}>m_{\ov{k}j}>m_{k+1,j+1}, 1\leq k\leq n-1, 1\leq j\leq k \}\cr
		      & &+ \#\{(k,j)\,|\, m_{\ov{k}j}>m_{kj}>m_{\ov{k},j+1}, 1\leq k\leq n, 1\leq j\leq k \}\,;\cr
	 R_o(GT)&=& \#\{(k,j)\,|\, m_{k+1,j}>m_{\ov{k}j}=m_{k+1,j+1}, 1\leq k\leq n-1, 1\leq j\leq k \}\,;\cr
	 L_e(GT)&=& \#\{(k,j)\,|\, m_{\ov{k}j}=m_{kj}>m_{\ov{k},j+1}, 1\leq k\leq n, 1\leq j\leq k-1 \}\,.\cr
\end{array}
\end{equation}
with 
\begin{equation}
      \x^{\xwgt(GT)} = \prod_{k=1}^n\ \prod_{j=1}^k x_k^{ 2m_{kj}-m_{\ov{k}j}-m_{\ov{k-1}j}} \,,
\end{equation}
and $m_{\ov{k},k+1}=0$ for $k=1,2,\ldots,n$.
\end{Corollary}

\noindent{\bf Proof}:
First it should be noted that
\begin{equation}
\begin{array}{l}
\ds \prod_{1\leq i\leq j \leq n}(x_i + q\,x_j)(1+q^{-1}x_i^{-1}x_j^{-1}) \cr
\ds = \prod_{i=1}^n (1+q)\,q^{-1}\,(qx_i+x_i^{-1}) \   \prod_{1\leq i<j \leq n} q^{-1} (x_i + q\,x_j)(q+ x_i^{-1}x_j^{-1}) \cr
\ds = \frac{(1+q)^n}{q^{n(n+1)/2}}\ \prod_{i=1}^n (q\,x_i+x_i^{-1}) \  \prod_{1\leq i<j \leq n} (x_i + q\,x_j)(q+ x_i^{-1}x_j^{-1})\,. \cr
\end{array}
\end{equation}
Second, that the $q$-dependence of $\wgt(GT)$ given by (\ref{eqn-q-dependence}) can be written in the form
\begin{equation} 
    \prod_{k=1}^n \left(\, (1+q) \prod_{j=1}^k q^{-1} \,\right) \times (1+q)^{B(GT)} \ q^{R_o(GT)+L_e(GT)}\,,
\end{equation}
where
\begin{equation}
\begin{array}{rcl}
        B(GT)&=& \sum_{k=1}^n \left(\,\sum_{j=1}^{k-1} (\chi(B_{kj})+\chi(B_{\ov{k}j})) +\chi(B_{kk})\,\chi(B_{\ov{k}k}) \right) \,:\cr
				R_o(GT)&=& \sum_{k=1}^n \sum_{j=1}^{k} \chi(R_{kj})\,;\cr
				L_e(GT)&=& \sum_{k=1}^n \sum_{j=1}^{k} \chi(L_{\ov{k}j})\,.\cr
\end{array}
\end{equation}
Writing these in the alternative forms displayed in (\ref{eqn-BRoLe}), noting that $R_o(GT)=$ the number of NW in odd rows of the CPM, and that $L_e(GT)=$ the number of NE in even rows of the CPM, and cancelling a common factor
\begin{equation}
      \prod_{k=1}^n \left(\, (1+q) \prod_{j=1}^k q^{-1} \,\right) =  \frac{(1+q)^n}{q^{n(n+1)/2}}
\end{equation}
that emerges on both sides of (\ref{eqn-Tok-GTq}) then yields (\ref{eqn-wgtGTq-final}). 
\qed
\medskip

Setting $y_k=qx_k$ for all $k$ in our previous example (\ref{eqn-example}) we find:
\begin{equation}
  B(GT)=\#\NS=7,\qquad R_o(GT)=\sum_{k=2}^5 \#\NW_{k} = 2, \qquad L_e(GT)=\sum_{k=1}^5 \#\NE_{\ov{k}} = 5
\end{equation}
so that if we exclude the overall factor of $(1+q)^n/q^{n(n+1)/2}=(1+q)^5/q^{15}$ we obtain
\begin{equation}\label{eqn-example-q}
  \qxwgt(GT) = (1+q)^{7}\, q^7\, x_2\,x_4^{-4} \,.
\end{equation}

\noindent{\bf Acknowledgements}
We wish to thank V. Gupta, U. Roy, and R. Van Peski for sending us their query and advanced notice of their results, including a statement of our Corollary~10, as well as further conjectures.

The first author (AMH) acknowledges the support of a
Discovery Grant from the Natural Sciences and Engineering Research Council of
Canada (NSERC). This work was supported by the Canadian Tri-Council Research Support Fund. 

\bigskip


\end{document}